\documentclass[12pt]{article}
\usepackage{amsmath}
\usepackage{amssymb}
\usepackage{dsfont}
\usepackage{cite}
\newsymbol \blackbox 1004
\newcommand{\eh}{\hfill}\newlength{\sperr}

\newenvironment{proof}{{\settowidth{\sperr}{\bf\rm
Proof}%
\par\addvspace{0.3cm}\noindent\parbox[t]{1.3\sperr}
{\bf\rm P\eh r\eh o\eh o\eh f\eh }%
}}{\nopagebreak\mbox{}
$\blackbox$\par\addvspace{0.3cm}}

\def\nn{\nonumber}

\def\vk{\varkappa}

\def\s{\sigma}
\def\la{\lambda}

\def\vt{\vartheta}

\def\wh{\widehat}
\def\wt{\widetilde}
\def\ov{\overline}

\def\p{\partial}

\def\BC{{\mathbb C}}

\def\BN{{\mathbb N}}

\def\cla{{\mathcal A}}

\def\cld{{\mathcal D}}

\def\clh{{\mathcal H}}

\def\cll{{\mathcal L}}

\def\clt{\mathcal{T}}

\def\cld{{\mathcal D}}

\newcommand{\E}{\mathrm{e}}
\newcommand{\I}{\mathrm{i}}

\newtheorem{Pa}{Paper}[section]
\newtheorem{Tm}[Pa]{{\bf Theorem}}

\newtheorem{Cy}[Pa]{{\bf Corollary}}
\newtheorem{Rk}[Pa]{{\bf Remark}}

\newtheorem{Pn}[Pa]{{\bf Proposition}}

\title{Dynamical canonical systems \\ and their explicit solutions}

\author{Alexander Sakhnovich}

\date{}
\begin{document}
\maketitle

\begin{abstract} Dynamical  canonical systems and their connections with
the classical (spectral) canonical systems are considered. We construct B\"acklund-Darboux transformation
and explicit solutions of the dynamical  canonical systems. We study also those properties
of the solutions, which are of interest in evolution and control theories.

\end{abstract}



\section{Introduction}
\setcounter{equation}{0}
In this paper we consider the system
\begin{align}& \label{D6}
 \frac{\p }{\p t} Y(x,t)=j\frac{\p }{\p x}\Big(\clh(x) Y(x,t)\Big), \quad \clh(x) >0, \quad j = \left[
\begin{array}{cc}
I_{m_1} & 0 \\ 0 & -I_{m_2}
\end{array}
\right],
\end{align}
where $I_p$ stands for the $p \times p$ idenity matrix, $m_1 \geq 0$, $m_2 \geq 0$, and we put $m:=m_1+m_2>0$, that is, $j$ is an $m\times m$ matrix
and $\clh(x)$ is an $m\times m$ matrix function.
The case $m_1>0$ and $m_2=0$ includes, for instance,
the transport equation and the case $m_1=m_2$ includes  the equation of a vibrating string (see, e.g., \cite{BMZ}).
We assume that $x$ takes values on a semiaxis
and without loss of generality consider the semiaxis $x \geq 0$.

In the case $m_1=m_2=p$,  we will show (using formal Fourier transformation)  an interesting correspondence between the system \eqref{D6} 
and the well-known canonical  system
\begin{align}& \label{D1}
\frac{d}{dx}y(x,\la)=\I \la J H(x)y(x,\la),   \quad H(x) \geq 0, \quad J=\left[
\begin{array}{cc}
0 & I_{p}  \\ I_{p} & 0
\end{array}
\right],
\end{align}
where  $\la \in \BC$ and $\BC$ stands for the complex plain.
Similar to the corresponding terminology for Dirac systems \cite{BeMi}, we call system \eqref{D6}  the {\it dynamical canonical system}
and system \eqref{D5} - the {\it spectral canonical system}. Since spectral Dirac (or Dirac-type) systems are equivalent to a subclass of 
spectral canonical systems (see, e.g., \cite[Sect. 1.1]{SaSaR}), dynamical Dirac systems  are equivalent to a subclass of 
dynamical canonical systems. (On the applications of the
dynamical Dirac systems see \cite{BeMi, BruK, Rak} and further references therein.)

Using B\"acklund-Darboux transformation for the  spectral canonical system \cite{SaA5} we construct B\"acklund-Darboux transformation 
and explicit solutions for the dynamical canonical system \eqref{D6}. (See Theorem \ref{TmDCGBDT} and Corollary \ref{ExplSol}, respectively.) 
We also study asymptotic behavior of these explicit solutions (see Theorem \ref{TmBeh}).

It is important that if  $j \clh_0(x)$ (where $\clh_0(x)$ is the Hamiltonian
of some initial dynamical system) is linear similar to some diagonal matrix $D(x)$, then the transformed  Hamiltonian $ \clh(x)$ has the same property,
that is, $j \clh(x)$ is similar 
to the same matrix $D(x)$. We note that the linear similarity of $j \clh$ to some diagonal matrix function $D$ is an essential assumption
in the main result (Theorem 1.5) in \cite{BMZ} on the operators $\big(\cla f\big)(x) =j\frac{d }{d x}  \clh(x) f(x)$ generating a $C_0$-semigroup
for \eqref{D6}. In fact, system \eqref{D7}  (which is more general than \eqref{D6}) is dealt with in \cite{BMZ}.

Each  generalized B\"acklund-Darboux transformation, in our GBDT version, is determined by some initial system and by some triple $\{A, S(0), \Pi(0)\}$ of parameter
matrices such that 
\begin{align}& \label{D8}
AS(0)-S(0)A^*=\I \Pi(0) j \Pi(0)^*.
\end{align}
We show that the energy of the constructed solutions is directly expressed via $A$ and $S(0)$ as well (see \eqref{D32}).

Further in the Introduction we discuss the mentioned above topics in greater detail and give various references.
We note that canonical (spectral canonical) systems
have been actively studied in the literature (see the books \cite{dB1, GoKr, SaL3, SaSaR}, quite recent papers \cite{Ach, Mogi, WiWo} and various references therein).
Using circumflex accent to denote
Fourier transform in the complex domain, that is, setting
\begin{align}& \label{D3}
\wh y(x,t)= \int_{-\infty}^{\infty}\E^{\I \la t}y(x, \la)d\xi \quad (\la = \xi + \I \eta),
\end{align}
and formally applying Fourier transform  to both sides of \eqref{D1}, we derive $\frac{\p }{\p x}\wh y(x,t)= J H(x)\frac{\p }{\p t}\wh y(x,t) $.
When the values of the matrix function (of Hamiltonian) $H(x)$ are strictly positive definite (i.e. $H(x)>0$), we set $Y(x,t)=H(x) \wh y(x,t)$ and rewrite
the last equation in an equivalent form
\begin{align}& \label{D5}
 \frac{\p }{\p t} Y(x,t)=J \frac{\p }{\p x}\Big(\clh(x) Y(x,t)\Big), \quad \clh(x) = H(x)^{-1}.
\end{align}

One may consider a more general case of \eqref{D1}, where $J$ is an arbitrary matrix such that  $J=J^*=J^{-1}$ and $J^*$ denotes the matrix adjoint to $J$.
Then the same formal Fourier transform of \eqref{D1} leads to system \eqref{D5} with the same (more general) $J$.
Since each matrix $J$, such that $J=J^*=J^{-1}$, is unitarily equivalent to some $j$ given in \eqref{D6},  system \eqref{D5}
is equivalent to the system \eqref{D6} with the corresponding $j$.

Transformations and explicit solutions of linear spectral systems and nonlinear evolution equations were actively studied in a number
of important works using commutation methods and various  versions of B\"acklund-Darboux transformations (see, e.g.,  \cite{D, Ge, GeT, MS, Gu, SaSaR}
and numerous references therein). Our GBDT version of the  B\"acklund-Darboux transformation was considered, for instance, in 
\cite{GKS6, KoSaTe, MST, SaA2, SaA5, SaA6, ALS10, SaSaR} (see also  further references therein). GBDT transformation for the
system \eqref{D1} (where $J=J^*=J^{-1}$) was introduced in \cite{SaA5}.

Using GBDT for  \eqref{D1}  and Fourier transformation \eqref{D3}, we could obtain explicit solutions of the dynamical system \eqref{D5}.
For example, in a related case of spectral and dynamical Dirac systems we established \cite{ALS-JMP} simple rigorous connections between Weyl functions of the spectral Dirac systems and
response functions of the dynamical Dirac systems   and used the procedure of solving inverse problem for the spectral Dirac system \cite{SaA7, ALS-JST}
in order to recover the dynamical Dirac system from its response function. However, in the case of canonical systems we  use GBDT for the system \eqref{D1} for heuristic
purposes only. In spite of similarity of many formulas in the dynamical and spectral cases, it seems simpler (as soon as one knows the result) to construct  GBDT and explicit solutions of the system \eqref{D6} directly. 

The condition that $\cla$ generates a $C_0$-semigroup (see \cite[Theorem 1.5]{BMZ}) which was discussed above, is closely connected
with the well-posedness of system \eqref{D6} with boundary control \cite{ZGMV}.
 We note that the study of well-posedness is basic in many evolution and control problems and a large number of recent
 publications is dedicated to this topic (see, e.g., \cite{BoKaR, ChuL, FeZ, LiR, MoP, RiCol, VaD}).
 
 GBDT and explicit solutions  for the dynamical canonical system are constructed in Section \ref{sec2}. 
 In Section 3 we
 study asymptotic behavior as well as energy and some other properties of the obtained solutions.
 We show also that $j \clh$ in our case is linear similar to simple diagonal matrices.
 The last section is Conclusion.

 \section{GBDT and explicit solutions \\ for the dynamical canonical system}\label{sec2}
\setcounter{equation}{0}
Further we fix an initial Hamiltonian $H\geq 0$, which (in the case $H>0$) is equivalent to fixing the initial Hamiltonian $\clh_0=H^{-1}$  (see \eqref{D5}).
Given an initial $m \times m$ Hamiltonian $H(x) \geq 0$ ($x\geq 0$), each GBDT is (as usual) determined by some $n \times n$ matrices $A$ and $S(0)>0$ ($n \in \BN$)
and by an $n \times m$ matrix $\Pi(0)$ which satisfy the matrix identity \eqref{D8},
where $j$ is defined in \eqref{D6}. Taking into account the initial values $\Pi(0)$ and $S(0)$ (and using the matrix $A$ and Hamiltonian $H(x)$) we introduce
matrix functions $\Pi(x)$ and $S(x)$ via the equations:
\begin{align}& \label{D9}
\Pi^{\prime}(x)=-\I A \Pi(x)jH(x) , \quad S^{\prime}(x)=\Pi(x)jH(x)j\Pi(x)^*,
\end{align}
where $\Pi^{\prime}:=\frac{d}{d x}\Pi$. It is easy to see \cite{SaA5} that \eqref{D8} and \eqref{D9} yield the identity
\begin{align}& \label{D10}
AS(x)-S(x)A^*\equiv \I \Pi(x) j \Pi(x)^*.
\end{align}
Since $S(0)>0$ and $S^{\prime}=\Pi jHj\Pi^* \geq 0$, the matrices $S(x)$ ($x\geq 0$) are positive definite and invertible. 
The so called Darboux matrix from Darboux transformations  is represented in GBDT (for each $\, x\, $)  as the transfer matrix function 
 in Lev Sakhnovich \cite{SaL1, SaL3, SaSaR} form:   $w_A(\la):=I_m-\I j \Pi^*S^{-1}(A-\la I_n)^{-1}\Pi$.
 More precisely, for the case of the canonical system (i.e., for the case that \eqref{D9} and \eqref{D10}
hold) we have \cite{SaA5}:
\begin{align}& \label{D11}
w_A^{\prime}(x,\la)=\big(\I \la j H(x)-\wt q_0(x)\big)w_A(x,\la)-\I \la w_A(x,\la) jH(x); 
\\ & \label{D12}
 \wt q_0(x):=j\Pi(x)^*S(x)^{-1}\Pi(x) j H(x)-j H(x) j \Pi(x)^*S(x)^{-1}\Pi(x),
\\ & \label{D13}
w_A(x,\la):=I_m-\I j \Pi(x)^*S(x)^{-1}(A-\la I_n)^{-1}\Pi(x).
\end{align}
\begin{Tm}\label{TmDCGBDT} Let a Hamiltonian $H(x)>0$  $(x\geq 0)$ and a triple of matrices
$A$, $S(0)>0$ and $\Pi(0)$, such that  \eqref{D8} holds, be given. Introduce matrix functions $\Pi(x)$, $S(x)$ and $u(x)$ via \eqref{D9}
and equalities
\begin{align}& \label{D14}
u^{\prime}(x)=-\wt q_0(x)u(x), \quad u(0)=I_m,
\end{align}
respectively, where $\wt q_0$ has the form \eqref{D12}. Then, the matrix functions
\begin{align}& \label{D15}
Y(x,t):=u(x)^*H(x)j\Pi(x)^*S(x)^{-1}\E^{\I tA}
\\ & \label{D16}
\clh(x):=u(x)^{-1}H(x)^{-1}\big((u(x)^*\big)^{-1}
\end{align}
satisfy dynamical canonical system \eqref{D6}.
\end{Tm}
\begin{proof}. It easily follows from \eqref{D12} that $\wt q_0^*j+j \wt q_0=0 $. Hence, formula \eqref{D14} implies that
$u^*ju\equiv j$, and so
\begin{align}& \label{D16'}
u(x)^{-1}=j u(x)^* j.
\end{align}
According to \eqref{D15} --\eqref{D16'} we have
\begin{align}& \label{D17}
\clh(x)Y(x,t)=u(x)^{-1}j\Pi(x)^*S(x)^{-1}\E^{\I tA}=ju(x)^*\Pi(x)^*S(x)^{-1}\E^{\I tA}.
\end{align}
Differentiating $\Pi^*S^{-1}$ and taking into account \eqref{D9} we obtain:
\begin{align}& \label{D18}
\big(\Pi^*S^{-1}\big)^{\prime}=\I H j \Pi^*A^*S^{-1}-\Pi^*S^{-1}\Pi j H j \Pi^* S^{-1}.
\end{align}
Multiplying both sides of \eqref{D10} by $S^{-1}$ from the left and from the right, we derive
\begin{align}& \label{D19}
A^*S^{-1}=S^{-1}A-\I S^{-1}\Pi j \Pi^*S^{-1}.
\end{align}
Next, we substitute \eqref{D19} into \eqref{D18}:
\begin{align}\nn
\big(\Pi^*S^{-1}\big)^{\prime}&=\I H j \Pi^*S^{-1}A+(H j \Pi^*S^{-1}\Pi j-\Pi^*S^{-1}\Pi j H j) \Pi^* S^{-1}
\\ & \label{D20}
=\I H j \Pi^*S^{-1}A+\wt q_0^* \Pi^* S^{-1}.
\end{align}
From \eqref{D14} and \eqref{D20} we see that $\big(u^*\Pi^*S^{-1}\big)^{\prime}
=\I u^* H j \Pi^*S^{-1}A$. Thus, formula \eqref{D17} yields
\begin{align}& \label{D21}
\frac{\p}{\p x}\left(\clh Y\right)= \I j u^*H j\Pi^*S^{-1}A \E^{\I tA}.
\end{align}
Finally, formula \eqref{D6} is immediate from  \eqref{D15} and \eqref{D21}.
\end{proof}
\begin{Rk}\label{Rklu}
 Comparing \eqref{D11} and \eqref{D14} we easily see that 
 \begin{align}& \label{D21!}
 u(x)=w_A(x,0)w_A(0,0)^{-1}.
 \end{align}
 We note that $w_A(x,0)$ and  $w_A(0,0)^{-1}$ are well-defined in the case $\det A\not= 0$.
Moreover, identity \eqref{D10} yields $w_A(x, \la)^{-1}=jw_A(x, \ov{\la})^*j$ $($see, e.g., \cite[f-la (13)]{SaA5}$)$. Thus, if only $\det A\not= 0$, the equality $u(x)=w_A(x,0)w_A(0,0)^{-1}$
and formula \eqref{D13} provide an explicit  expression for $u$ in terms of $S$ and $\Pi \, :$
\begin{align}& \label{D22}
u(x)=\big(I_m-\I j \Pi(x)^*S(x)^{-1}A^{-1}\Pi(x)\big)\big(I_m +\I j \Pi(0)^*\big(A^{-1}\big)^*S(0)^{-1}\Pi(0)\big).
\end{align}
\end{Rk}

Let us consider in detail the case of the {\it trivial} initial Hamiltonian, that is, the case $H(x) \equiv I_m$. In order to avoid confusion, we assume for this case that $m_1>0$ and $m_2>0$ in \eqref{D6}.
(The modification of our considerations for $m_1>0$ and $m_2=0$ is simple and evident.)
Partition $\Pi(0)$ into two blocks $\Pi(0)=\begin{bmatrix} \vt_1 & \vt_2 \end{bmatrix}$
(where $\vt_i$ are $n \times m_i$ matrices) and derive from \eqref{D9} explicit expressions for $\Pi(x)$ and $S(x)$:
\begin{align}& \label{D23}
\Pi(x)=\begin{bmatrix}\E^{-\I xA} \vt_1 & \E^{\I xA} \vt_2 \end{bmatrix}, \quad S(x)=S(0)+\int_0^x\Pi(\xi)\Pi(\xi)^*d\xi.
\end{align}
The following corollary of Theorem \ref{TmDCGBDT} presents explicit solutions of \eqref{D6}.
\begin{Cy}\label{ExplSol} Let  a triple of matrices
$\{A, \, S(0), \,\Pi(0)\}$, such that $S(0)>0$, $\det A \not=0$ and \eqref{D8} holds, be given. Introduce matrix functions $\Pi(x)$, $S(x)$ and $u(x)$ via 
\eqref{D23} and \eqref{D22}, respectively. Then, the matrix functions
\begin{align}& \label{D24}
Y(x,t)=u(x)^*j\Pi(x)^*S(x)^{-1}\E^{\I tA},
\quad
\clh(x)=u(x)^{-1}\big((u(x)^*\big)^{-1}
\end{align}
satisfy dynamical canonical system \eqref{D6}.
\end{Cy}
In the formula \eqref{D23}, the expression for $\Pi$ is somewhat more explicit then the expression for $S$.
However, the expression for $S$ in \eqref{D23} is easily rewritten in terms of the solutions $C_k$ of the
matrix identities (linear algebraic equations):
\begin{align}& \label{D23'}
AC_k-C_kA^*=\I \vt_k \vt_k^* \quad (k=1,2).
\end{align}
Indeed, in view of \eqref{D23'} and the first equality in \eqref{D23} we have
\begin{align}& \label{D23+}
\big(\E^{-\I xA} C_1\E^{\I xA^*}-\E^{\I xA} C_2\E^{-\I xA^*}\big)^{\prime}=\Pi(x)\Pi(x)^*.
\end{align}
Using \eqref{D23+}, we immediately simplify the expression for $S$ in \eqref{D23} and derive the corollary below.
\begin{Cy} \label{CyExpS} Let \eqref{D23'} hold. Then $S(x)$ given by the second equality in \eqref{D23} takes the form
\begin{align}& \label{D23!}
S(x)=S(0)-C_1+C_2+\E^{-\I xA} C_1\E^{\I xA^*}-\E^{\I xA} C_2\E^{-\I xA^*}.
\end{align}
\end{Cy}
\begin{Rk}\label{RkExpl} 
Clearly, \eqref{D8} implies that, if \eqref{D23'} holds for $k=1$, we may set $C_2=C_1-S(0)$,
and \eqref{D23'} will hold for $k=2$ $($with this choice of $C_2)$ as well. If \eqref{D23'} holds for $k=2$ we may set $C_1=C_2+S(0)$.

We note that the identities \eqref{D23'}  always have unique solutions $C_k$ when $\s(A)\cap \s(A^*)=\emptyset$.

Substituting  the first equality in \eqref{D23} and equality \eqref{D23!}  $($i.e., explicit expressions for $\Pi(x)$ and $S(x))$ into \eqref{D22},
we obtain an explicit expression for $u(x)$. Now, substituting explicit expressions for $\Pi(x)$, $S(x)$ and $u(x)$ into \eqref{D24}, we obtain
explicit expressions for $\clh(x)$ and $Y(x,t)$.
 \end{Rk}

Explicit solutions of the related to the system \eqref{D6} Loewner's system 
$$\frac{\p}{\p x}Y (x,t)=\cll (x,t) \frac{\p}{\p t}Y (x,t)$$ 
were
constructed in \cite{FKRS}. However, the dependence of $\cll $ on both
variables is essential in that construction and the procedure itself is more complicated.


\section{Analysis of the obtained  solutions}\label{sec3}
\setcounter{equation}{0}
\paragraph{1.} In this section, we consider $\clh$ and $Y$ constructed in Theorem \ref{TmDCGBDT}.
According to \cite[Section 2]{ZGMV}, the energy $E_{Yh}(t)$ of the solution $Y(x,t)h$ ($h\in \BC^m$) on the
interval $0\leq x \leq a$ is  given in physical models by the expression
\begin{align}& \label{D29}
E_{Yh}(t)^2=\int_0^a h^*Y(x,t)^*\clh(x)Y(x,t)hdx.
\end{align}
The vector $h$ in \eqref{D29} determines the boundary condition for the solution $Yh$:
\begin{align}& \label{D29'}
\big(Yh\big)(0,t)=H(0)j\Pi(0)^*S(0)^{-1}\E^{\I tA}h.
\end{align}

In view of \eqref{D15} and \eqref{D16} we have
\begin{align}& \label{D30}
 Y(x,t)^*\clh(x)Y(x,t)=\E^{-\I tA^*}S(x)^{-1}\Pi(x)jH(x)j\Pi(x)^*S(x)^{-1}\E^{\I tA}.
\end{align}
The second relation in \eqref{D9} implies that $S^{-1}\Pi jH j\Pi^*S^{-1}=-\big(S^{-1}\big)^{\prime}$.
Hence, we rewrite \eqref{D30} in the form
\begin{align}& \label{D31}
 Y(x,t)^*\clh(x)Y(x,t)=-\E^{-\I tA^*}\big(S(x)^{-1}\big)^{\prime}\E^{\I tA}.
\end{align}
Substituting \eqref{D31} into \eqref{D29} we obtain our next statement.
\begin{Pn}\label{PnE} Let the conditions of Theorem \ref{TmDCGBDT} hold. Then the energy $E_{Yh}$, where
$Y$ is constructed in Theorem \ref{TmDCGBDT}, is given by the formula
\begin{align}& \label{D32}
E_{Yh}(t)=\sqrt{h^*\E^{-\I tA^*}\big(S(0)^{-1}-S(a)^{-1}\big)\E^{\I tA}h}.
\end{align}
\end{Pn}
We note that the energy in \eqref{D32} is directly expressed via the boundary values $S(0)$ and $S(a)$.

Let us set
\begin{align} & 
& \label{D33}
s(x, Y, h)=h^*Y(x,t)^*u(x)^{-1}H(x)^{-1}jH(x)^{-1}\big(u(x)^*\big)^{-1}Y(x,t)h.
\end{align}
Then, relation \eqref{D15} 
implies that
\begin{align} & \label{D34}
s(x,Y,h)=h^*\E^{-\I tA^*}S(x)^{-1}\Pi(x)j\Pi(x)^*S(x)^{-1}\E^{\I tA}h.
\end{align}
On the other hand, \eqref{D19} yields the formula 
\begin{align}\nn
\frac{\p}{\p t}\left(h^*\E^{-\I tA^*}S(x)^{-1}\E^{\I tA}h\right)&=-\I h^*\E^{-\I tA^*}(A^*S(x)^{-1}- S(x)^{-1}A)\E^{\I tA}h
\\ & \label{D35}
=-h^*\E^{-\I tA^*}S(x)^{-1}\Pi(x)j\Pi(x)^*S(x)^{-1}\E^{\I tA}h.
\end{align}
Formulas \eqref{D32}--\eqref{D35} imply the following equality, which is similar
to the relations for supply rate for port-Hamiltonian systems that appear in the literature:
\begin{align}
 & \label{D36}E_{Yh}(t_2)^2- E_{Yh}(t_1)^2 =
\int_{t_1}^{t_2}\Big(s(a,Y,h)-s(0,Y,h)\Big)dt.
\end{align}
Recall that $S(x)>0$ and that $S(x)$ is a monotonically nondecreasing matrix function. 
Hence, there is a finite nonnegative limit $\vk_S$ of $S(x)^{-1}$ when $x $ tends to infiity.
\begin{Rk}\label{RkAsS} When we consider the case of the semiaxis $0 \leq x <\infty$
instead of the interval $0 \leq x \leq a$, we substitute $\vk_S$ $($instead of $S(a)^{-1})$
into the expression \eqref{D32} for energy. 
\end{Rk}

\paragraph{2.} 
In this paragraph, we study the behavior of the explicit solution $Y(x,t)$ given by \eqref{D24}.
According to Corollary \ref{ExplSol}, $Y$ is expressed in terms of $\Pi(x)$ and $S(x)$.
Similar to \cite[Sect. 6]{GKS6} (see also \cite{SaA8}) we include into consideration
matrix functions
\begin{align}
 & \label{A1} 
 Q(x)=\E^{\I x A}S(x)\E^{-\I x A^*}, \quad R(x)=\E^{-\I x A}S(x)\E^{\I x A^*} \qquad (x\geq 0).
 \end{align}
 Taking into account \eqref{D10}, \eqref{D23} and \eqref{A1} we have
\begin{align}
\nn
 Q^{\prime}(x)&=\E^{\I x A}\big(S^{\prime}(x)+\I(AS(x)-S(x)A^*)\big)\E^{-\I x A^*}
\\ &
  \label{A2}  
 =\E^{\I x A}\Pi(x)(I_m-j)\Pi(x)^*\E^{-\I x A^*}\geq 0,
\\
\nn
 R^{\prime}(x)&=\E^{-\I x A}\big(S^{\prime}(x)-\I(AS(x)-S(x)A^*)\big)\E^{\I x A^*}
 \\ &
  \label{A3}  
 =\E^{-\I x A}\Pi(x)(I_m+j)\Pi(x)^*\E^{\I x A^*} \geq 0.
  \end{align}
  It is immediate from \eqref{A1}--\eqref{A3} that $Q$ and $R$ are positive-definite and monotonically nondecreasing matrix functions,
  and so the limits $Q(x)^{-1}$ and $R(x)^{-1}$ ($x\to \infty$) exist. We set
  \begin{align}
 & \label{A4} 
\vk_Q:=\lim_{x\to \infty} Q(x)^{-1}, \quad \vk_R:=\lim_{x\to \infty} R(x)^{-1}.
 \end{align}
 Relations \eqref{D13}, \eqref{A1} and the first equality in \eqref{D23} imply that $w_A(x,0)$ admits representation
  \begin{align}
 & \label{A5} 
w_A(x,0)=\begin{bmatrix} I_{m_1}-\I\vt_1^*Q(x)^{-1}A^{-1}\vt_1 &-\I \vt_1^*\E^{2\I x A^*}R(x)^{-1}A^{-1}\vt_2 \\
\I \vt_2^*\E^{-2\I x A^*}Q(x)^{-1}A^{-1}\vt_1 &  I_{m_2}+\I\vt_2^*R(x)^{-1}A^{-1}\vt_2 
\end{bmatrix}.
 \end{align}
 In view of \eqref{D21!}, the asymptotics of  $w_A(x,0)$, which is given in the next proposition, provides
 the asymptotics of $u$ in the expression \eqref{D24} for $Y$.
 \begin{Pn}\label{Pnwa} Let  a triple of matrices
$\{A, \, S(0), \,\Pi(0)\}$, such that \eqref{D8} holds, $S(0)>0$ and $\det A \not=0$, be given. Introduce matrix functions $\Pi(x)$, $S(x)$ and $w_A(x,\la)$ via 
\eqref{D23} and \eqref{D13}, respectively. Then we have
   \begin{align}
 & \label{A6} 
\lim_{x\to \infty} w_A(x,0)=\begin{bmatrix} I_{m_1}-\I\vt_1^*\vk_Q A^{-1}\vt_1 & 0 \\
0 &  I_{m_2}+\I\vt_2^* \vk_R A^{-1}\vt_2 
\end{bmatrix}.
 \end{align}
 \end{Pn}
 \begin{proof}.  The equality \eqref{A6} is immediate from \eqref{A4}, \eqref{A5} and relations
  \begin{align}& \label{A7} 
\lim_{x\to \infty}\vt_2^*\E^{-2\I x A^*}Q(x)^{-1}=0, \quad \lim_{x\to \infty} \vt_1^*\E^{2\I x A^*}R(x)^{-1}=0,
\end{align} 
which we prove below. We prove  \eqref{A7} similar to \cite[Prop. 3.1]{GKS6}, although relations \eqref{A7}
are derived here in a somewhat more general situation. First note that formulas \eqref{D10} and \eqref{A1} imply that
  \begin{align}& \label{A8} 
AQ(x)-Q(x)A^*=\I \E^{\I x A}\Pi(x)j\Pi(x)^*\E^{-\I x A^*} , \\ 
& \label{A9}
AR(x)-R(x)A^*=\I \E^{-\I x A}\Pi(x)j\Pi(x)^*\E^{\I x A^*}.
\end{align} 
Multiplying \eqref{A8} by $Q(x)^{-1}$ and \eqref{A9} by $R(x)^{-1}$ from both sides and using the expression
for $\Pi$ from \eqref{D23}, we obtain
  \begin{align}& \label{A10} 
Q^{-1}A-A^*Q^{-1}-\I Q^{-1}\vt_1\vt_1^* Q^{-1} =-\I Q^{-1} \E^{2\I x A}\vt_2 \vt_2^*\E^{-2 \I x A^*}Q^{-1} , \\ 
& \label{A11}
R^{-1}A-A^*R^{-1}+\I R^{-1}\vt_2\vt_2^* R^{-1} =\I R^{-1} \E^{-2\I x A}\vt_1 \vt_1^*\E^{2 \I x A^*}R^{-1}.
\end{align}
On the other hand,  \eqref{A2}, \eqref{A3} and the first equality in \eqref{D23} yield
  \begin{align}& \label{A12} 
\big(Q^{-1}\big)^{\prime} (x)=-2 Q(x)^{-1}  \E^{2\I x A}\vt_2 \vt_2^*\E^{-2 \I x A^*}        Q(x)^{-1},
\\ & \label{A13}
\big(R^{-1}\big)^{\prime} (x)=-2 R(x)^{-1}  \E^{-2\I x A}\vt_1 \vt_1^*\E^{2 \I x A^*}        R(x)^{-1}.
\end{align}
Equality \eqref{A12} and the first equality in \eqref{A4} imply that the entries of $\big(Q^{-1}\big)^{\prime}$ are summable
on $[0,\, \infty)$ (i.e., $\big(Q^{-1}\big)^{\prime}$ is summable on $[0,\, \infty)$). Hence, taking into account that $Q^{-1}$ is bounded,
we see that the derivative of the left-hand side of \eqref{A10} is summable, and so the derivative of the right-hand side of \eqref{A10} is summable.

Moreover, since $\big(Q^{-1}\big)^{\prime}$ is summable, formula \eqref{A12} implies that the right-hand side of \eqref{A10} is summable.

The fact that the  right-hand side of \eqref{A10} and its derivative are both summable on $[0,\, \infty)$ means that the right-hand side of \eqref{A10}
tends to zero, that is,
  \begin{align}& \label{A14} 
\lim_{x\to \infty}Q^{-1} \E^{2\I x A}\vt_2 \vt_2^*\E^{-2 \I x A^*}Q^{-1}=0.
\end{align} 
The first equality in \eqref{A7} follows from \eqref{A14}. The second equality in \eqref{A7}  is proved in the same way using \eqref{A11} and
\eqref{A13}.
 \end{proof}
Taking into account  \eqref{D24} and the proposition above we show that the behavior of $Y(x,t)$ is characterized by two exponents
$\E^{\I(t+x)A}$ and $\E^{-\I(t+x)A}$.
\begin{Tm}\label{TmBeh} Let the conditions of Corollary \ref{ExplSol} hold. Then we have
  \begin{align}\nn
Y(x,t)=\bigg(& jw_A(0,0)\begin{bmatrix} I_{m_1} 
+\I\vt_1^*\big(A^{-1}\big)^*\vk_Q \vt_1 & 0 \\
0 &  I_{m_2}-\I\vt_2^*\big(A^{-1}\big)^* \vk_R \vt_2 
\end{bmatrix} 
\\ & \label{A15}
+o(1)\bigg)
\begin{bmatrix}(\vt_1^*\vk_Q+o(1))\E^{\I(t+x)A} \\
(\vt_2^*\vk_R+o(1))\E^{\I(t-x)A}
\end{bmatrix} ,
\end{align} 
where the  functions  ``little-o" tend to $0$ when $x$ tends to infinity.
\end{Tm}
\begin{proof}. According to Proposition \ref{Pnwa} and Remark \ref{Rklu}, the first factor (in round brackets) on the right-hand side
of \eqref{A15} gives us $u(x)^*j$. 

Using the first equality in \eqref{D23} and relations \eqref{A1} and \eqref{A4}, we have
\begin{align}& \label{A16}
\Pi(x)^*S(x)^{-1}=\begin{bmatrix}\vt_1^*Q(x)^{-1}\E^{\I xA} \\
\vt_2^*R(x)^{-1}\E^{-\I xA}
\end{bmatrix}=\begin{bmatrix}\big(\vt_1^*\vk_Q+o(1)\big)\E^{\I xA} \\
\big(\vt_2^*\vk_R +o(1)\big)\E^{-\I xA}
\end{bmatrix}.
\end{align}
Hence, we rewrite $\Pi(x)^*S(x)^{-1}\E^{\I t A}$ in the form of the second factor on the right-hand side
of \eqref{A15}. Now, \eqref{A15} is immediate from \eqref{D24}.
\end{proof}
\begin{Rk} Since $R(x)^{-1}=\E^{-2\I xA^*}Q(x)^{-1}\E^{2\I xA}$ and \eqref{A7} is valid, we can write a representation
of $\Pi^*S^{-1}$ in a form somewhat different from \eqref{A16}. Namely, we obtain$:$
\begin{align}& \label{A17}
\Pi(x)^*S(x)^{-1}=\begin{bmatrix}\vt_1^*Q(x)^{-1}\\
\vt_2^*\E^{-2\I xA^*}Q(x)^{-1}
\end{bmatrix}\E^{\I xA}=\begin{bmatrix}\big(\vt_1^*\vk_Q+o(1)\big) \\
o(1)
\end{bmatrix}\E^{\I xA}.
\end{align}
Correspondingly, we rewrite \eqref{A15} in the form, which is convinient for the study of asymptotics
of $Y$ with large $x$ and $t$, when the spectrum of $A$ belongs to the upper semiplane$:$
  \begin{align}\nn
Y(x,t)=\bigg(& jw_A(0,0)\begin{bmatrix} I_{m_1} 
+\I\vt_1^*\big(A^{-1}\big)^*\vk_Q \vt_1 & 0 \\
0 &  I_{m_2}-\I\vt_2^*\big(A^{-1}\big)^* \vk_R \vt_2 
\end{bmatrix} 
\\ & \label{A18}
+o(1)\bigg)
\begin{bmatrix}\big(\vt_1^*\vk_Q+o(1)\big) \\
o(1)
\end{bmatrix}\E^{\I (x+t)A}.
\end{align} 
\end{Rk}
Clearly, we can write an analog of \eqref{A17} with $R$ instead of $Q$ as well.
\paragraph{3.} 
We mentioned already in the Introduction that the linear similarity of $j \clh(x)$ to some diagonal matrix function $\cld(x)$ and
the corresponding representation  $j \clh(x)=\clt(x)\cld(x)\clt(x)^{-1}$ is essential for well-posedness problems.
 In our case, using \eqref{D16'} we can rewrite \eqref{D16} in the form
\begin{align}& \label{D25}
j \clh(x)=\big(j u(x)j)^{-1}j H(x)^{-1}ju(x)j.
\end{align}
Hence, if $jH^{-1}$ admits a representation
\begin{align}& \label{D26}
j H(x)^{-1}=T(x)D(x)T(x)^{-1},
\end{align}
where $D(x)$  is some diagonal matrix function, equalities \eqref{D25} and \eqref{D26} yield
the following result.
\begin{Tm}\label{TmSim} 
Let the conditions of Theorem \ref{TmDCGBDT} and equality \eqref{D26} hold. Then we have
\begin{align}& \label{D27}
j \clh(x)= \clt(x)D(x)\clt(x)^{-1} \quad {\mathrm{for}} \quad \clt(x)=\big(j u(x)j)^{-1}T(x).
\end{align}
\end{Tm}
\begin{Rk}\label{RkSim}
In the case of the explicit solutions considered in Corollary \ref{ExplSol} we have $H(x)=I_m$ and so
\begin{align}& \label{D28}
j \clh(x)= \clt(x)j \clt(x)^{-1} \quad {\mathrm{for}} \quad \clt(x)=\big(j u(x)j)^{-1}.
\end{align}
\end{Rk}
In the case treated in Remark \ref{RkSim} we  obtain explicit expressions for
the spans $Z^{+}(x)$ and $Z^{-}(x)$ of the eigenvectors of $j\clh(x)$ corresponding to the positive and negative, respectively,  eigenvalues of $j\clh(x)$:
\begin{align}&  \label{Z+}
Z^{+}(x)=(ju(x)j)^{-1}(I_m+j)\BC^m=u(x)^*(I_m+j)\BC^m, 
\\ &   \label{Z-}
 Z^{-}(x)=u(x)^*(I_m-j)\BC^m.
\end{align}
Using \eqref{Z+} and \eqref{Z-} one can easily check the validity of the condition~2. from \cite[Theorem 1.5]{BMZ} that the operator $j\frac{d}{d x}\clh(x)$ is a generator
of  a $C_0$-semigroup.
\section{Conclusion}\label{concl}
\setcounter{equation}{0}
The  solutions of the dynamical canonical system, which are constructed in this note, and their behavior
are of independent interest. These solutions may also serve as examples
 for some problems arising in control theory.

The results of the note are easily generalized
for the case  $\clh = \clh^*$ where the requirement $\clh >0$ is omitted.
We could also study (in the spirit of \cite{SaA8}) blow up solutions with singularities, which appear, if we omit
the requirement $S(0)>0$. We mention that GBDT was successfully applied for the
construction of explicit solutions of nonlinear dynamical systems as well (see, e.g., various references in \cite{FKRS, SaA6, ALS10, SaSaR}).

In the control problems for \eqref{D6} boundary conditions are essential.
\begin{Rk} \label{RkBC} 
Boundary conditions of the form 
\begin{align}& \label{BC1}
\wt W\begin{bmatrix} \clh(a,t)Y(a,t) \\  \clh(0,t)Y(0,t)\end{bmatrix}=0,
\end{align}
where $\wt W$ is an $m \times 2m$ matrix, are considered in \cite{BMZ}. Choosing $m$-dimensional
invariant subspaces $L_A$ of $A$ and using \eqref{D17}, we can always construct such matrices
$\wt W$ that \eqref{BC1} holds for all  the solutions $Y(x,t)h$ $(h \in L_A)$, where
$Y$ is given by \eqref{D15}.
\end{Rk} 

It is our plan  to consider (in the next paper) B\"acklund-Darboux transformations for the case of a more general than \eqref{D6} system
\begin{align}& \label{D7}
 \frac{\p }{\p t} Y(x,t)=\Big(P_1\frac{\p }{\p x} +P_0)\Big(\clh(x) Y(x,t)\Big), \quad \clh(x) >0, \quad P_1=P_1^*,
 \end{align}
and study applications to control theory in a more detailed way. 

\bigskip

\noindent{\bf Acknowledgments.}
 {This research   was supported by the
Austrian Science Fund (FWF) under Grant  No. P24301.}
\newpage

\begin{flushright}

A.L. Sakhnovich,\\
Institute for Analysis and Scientific Computing,\\
Vienna
University
of
Technology, \\
Wiedner Haupstr. 8-10 / 101, 1040 Wien,
Austria, \\
e-mail: {\tt oleksandr.sakhnovych@tuwien.ac.at}

\end{flushright}


\end{document}